\documentclass[12pt]{amsart}
\author{Kai (Steve) Fan}
\author{Paul Pollack}
\address{Department of Mathematics\\ University of Georgia\\ Athens, GA 30602}
\email{Steve.Fan@uga.edu}
\email{pollack@uga.edu}
\title{Extremal elasticity of quadratic orders}
\usepackage{amsmath,amssymb,amsthm,microtype,geometry,mathscinet,enumitem,url}
\usepackage{parskip}
\usepackage[utf8]{inputenc}
\usepackage{mathtools}

\renewcommand\subset\subseteq
\newcommand{\verteq}{\rotatebox{90}{$\,=$}}
\newcommand{\vertcong}{\rotatebox{-90}{$\cong\,$}}

\newcommand{\equalto}[2]{\underset{\scriptstyle\overset{\mkern4mu\verteq}{#2}}{#1}}
\newcommand{\conggto}[2]{\underset{\scriptstyle\overset{\mkern4mu\vertcong}{#2}}{#1}}

\DeclareMathAlphabet{\curly}{U}{rsfs}{m}{n}
\newtheorem{thm}{Theorem}[section]

\newtheorem{prop}[thm]{Proposition}
\newtheorem{lem}[thm]{Lemma}

\newtheorem*{thmA}{Theorem A}
\newtheorem*{thmB}{Theorem B}
\newtheorem*{thmC}{Theorem C}
\newtheorem*{thmD}{Theorem D}
\newtheorem*{thmE}{Theorem E}

\theoremstyle{remark}
\newtheorem*{rmk}{Remark}

\DeclareMathOperator{\rk}{\mathrm{rk}}
\DeclareMathOperator{\Dav}{\mathrm{Dav}}
\DeclareMathOperator{\Rad}{\mathrm{Rad}}
\DeclareMathOperator{\Exp}{\mathrm{Exp}}
\DeclareMathOperator{\lcm}{\mathrm{lcm}}
\def\colonequals{:=}
\def\sgn{\mathrm{sgn}}
\def\Dd{\mathcal{D}}

\begin{document}
\renewcommand{\labelenumi}{(\roman{enumi})}
\newcommand{\genlegendre}[4]{%
  \genfrac{(}{)}{}{#1}{#3}{#4}%
  \if\relax\detokenize{#2}\relax\else_{\!#2}\fi
}
\newcommand{\leg}[3][]{\genlegendre{}{#1}{#2}{#3}}
\newcommand{\dlegendre}[3][]{\genlegendre{0}{#1}{#2}{#3}}
\newcommand{\tlegendre}[3][]{\genlegendre{1}{#1}{#2}{#3}}\def\Ll{\mathcal{L}}
\newcommand{\psmod}[1]{\,(\textup{\text{mod}}\,{#1})}
\def\N{\mathbb{N}}
\def\Aa{\mathcal{A}}
\def\I{\mathcal{I}}
\def\Q{\mathbb{Q}}
\def\Oo{\mathcal{O}}
\def\Nm{\mathrm{Nm}}
\newcommand\rad{\mathrm{rad}}
\def\Z{\mathbb{Z}}
\def\F{\mathbb{F}}
\def\R{\mathbb{R}}
\def\C{\mathbb{C}}
\def\Cl{\mathrm{Cl}}
\def\PrinCl{\mathrm{PrinCl}}
\def\PreCl{\mathrm{PreCl}}
\def\Pp{\mathcal{P}}
\def\Ss{\mathcal{S}}
\def\Uu{\mathcal{U}}
\def\th{\mathrm{th}}
\def\Frob{\mathrm{Frob}}
\newcommand\Li{\mathrm{Li}}
\begin{abstract} We study how large and small elasticity can be for orders belonging to a fixed quadratic field, in terms of the corresponding conductors. For example, we show that if $K$ is an imaginary quadratic field, then the order of conductor $f$ in $K$ has elasticity exceeding $(\log{f})^{c_1 \log\log\log{f}}$ for all $f$ that are sufficiently large. On the other hand, this elasticity is smaller than $(\log{f})^{c_2\log\log\log{f}}$ for infinitely many $f$. Here $c_1, c_2$ are universal positive constants. The proofs borrow methods from analytic number theory previously employed to study statistics of the multiplicative groups $(\Z/m\Z)^{\times}$.
\end{abstract}

\maketitle

\section{Introduction}
Let $R$ be an \textsf{atomic domain}, i.e., an integral domain in which every nonzero nonunit factors as a product of irreducibles. The \textsf{elasticity} of $R$, denoted $\rho(R)$, is the supremum of all ratios $r/s$, where $r$ and $s$ are positive integers for which there exists an equation
\[ \pi_1 \cdots \pi_r = \rho_1 \cdots \rho_s \]
with each $\pi_i$ and $\rho_j$ irreducible in $R$. So for instance, every unique factorization domain (not a field) has elasticity $1$. Conversely, if $R$ is an atomic domain with elasticity $1$, then $R$ might be thought of as `halfway' to unique factorization: any two factorizations of the same nonzero nonunit involve the same number of irreducibles, although the irreducibles themselves need not be pairwise associate. A domain of elasticity $1$ is called a \textsf{half-factorial domain} (or \textsf{HFD}).

The study of half-factorial domains was initiated in a 1960 paper of Carlitz \cite{carlitz60}, where it is shown that the ring of integers $\Oo_K$ of the number field $K$ is a half-factorial domain precisely when the corresponding class number $h_K = 1$ or $2$. The general notion of elasticity was introduced by Valenza in a 1990 paper\footnote{received by the journal in 1980 !} \cite{valenza90} and further studied by Narkiewicz \cite{narkiewicz95} and Steffan \cite{steffan86}. Collecting the results of these three papers, one has a complete determination of the elasticity of the ring of integers of an arbitrary number field.

For a finite abelian group $G$, we write $\Dav G$ for the \textsf{Davenport constant} of $G$, meaning the least positive integer $D$ with the following property: Any sequence $g_1, \dots, g_D$ of elements of $G$ possesses a nonempty subsequence multiplying to the identity. 
\begin{thmA}[Narkiewicz--Steffan--Valenza] Let $R$ be a Dedekind domain with finite class group. Suppose that every ideal class of $R$ is represented by at least one maximal ideal of $R$. Then
\[ \rho(R) = \max\left\{1, \frac{1}{2}\Dav \Cl(R)\right\}, \]
where $\Cl(R)$ is the class group of $R$.
\end{thmA}
\noindent Results of class field theory guarantee that the hypotheses of Theorem A are satisfied for any ring of $S$-integers in any global field. Thus, Theorem A gives a complete description of elasticities for the rings most prominent in number theory. 

Of course, there are many other important rings in algebra and number theory, and there is a rich literature investigating elasticity and half-factoriality more generally. See \cite{anderson97} and \cite{cc00} for surveys. In this paper, which is a companion piece to our recent article \cite{FP}, we will be examining elasticity in orders of quadratic fields.

Let $K$ be a quadratic field (a field extension of $\Q$ with $[K:\Q]=2$). An \textsf{order} in $K$ is a subring of $\Oo_K$ strictly larger than $\Z$. The orders in $K$ are in one-to-one correspondence with the natural numbers (positive integers) $f$. Each order $\Oo$ in $K$ has the form
\[ \Oo_f := \{\alpha \in \Oo_K: \alpha\equiv a\bmod{f\Oo_K}\text{ for some rational integer $a$} \} \]
for a unique $f\in \N$; conversely, $\Oo_f$ as defined above is always an order in $K$. It sometimes helpful to have a more explicit description of $\Oo_f$: It is well-known that for each quadratic field $K$, there is a unique squarefree integer $D$ with $K= \Q(\sqrt{D})$. If we set
\begin{equation}\label{eq:taudef} \tau_D = \begin{cases} \sqrt{D} &\text{if $D\equiv 2,3\pmod{4}$}, \\
\frac{1}{2} (1+\sqrt{D}) &\text{if $D\equiv 1\pmod{4}$},
    \end{cases}
\end{equation}
then $\Oo_f = \Z + f\tau_D \Z$. We refer to $\Oo_f$ as the order of \textsf{conductor} $f$ in $K$. Note that if $f\mid f'$, then $\Oo_{f'} \subset \Oo_{f}$, and that when $f=1$, we have $\Oo_1 = \Oo_K$; we will refer to $\Oo_K$ as the \textsf{maximal order} in $K$. For more on the basic theory of quadratic orders, see the lovely  book of Cox \cite{cox22} (p.\ 105 and following).

It was noted by Zaks \cite{zaks76, zaks80} that  $\Z[\sqrt{-3}]$, the order of conductor $2$ in $K=\Q(\sqrt{-3})$, is a half-factorial domain. However, half-factorial quadratic orders were not systematically studied until the later work of Coykendall \cite{coykendall01} and Halter-Koch \cite{HK83}. (This latter paper of Halter-Koch is not easy to come by; a more readily-available source for this material is \cite{GHK06}; see p.\ 226 and following.) These papers contain a complete algebraic characterization of half-factorial quadratic orders.

The characterization in the imaginary case is surprisingly simple: $\Z[\sqrt{-3}]$ is the unique nonmaximal half-factorial quadratic order \cite[Theorem 2.3]{coykendall01}. (Recall that maximal orders fall under the purview of Carlitz's 1960 theorem, so we may always restrict to the nonmaximal case.) The real case is substantially more complicated. Here there are various ways to state the characterization, some more explicit than others; the version we quote below is based on Coykendall's paper \cite{coykendall01} and a subsequent manuscript of Coykendall--Malcolmson--Okoh \cite{CMO17}.

\begin{thmB} Let $\Oo$ be the order of conductor $f$ in the real quadratic field $K$. In order for $\Oo$ to be half-factorial, it is necessary that
\begin{enumerate}[label=\normalfont(\roman*)]
\item $\Oo_K$ be half-factorial {\rm (}equivalently, $h_K=1$ or $h_K=2${\rm )}, \emph{and}
\item $f=p$ or $f=2p$, where $p$ is a prime inert in $K$; if $f=2p$, we require $p > 2$ and both $2$ and $p$ to be inert in $K$.
\end{enumerate}
Conversely, suppose $\Oo_K$ is half-factorial and $f=p$ is a prime inert in $K$. Let $\epsilon$ be the fundamental unit of $K$.\footnote{While not important here, we mention for later: We always view quadratic fields as embedded in $\C$. If $K=\Q(\sqrt{D})$ is real-quadratic, we assume $\sqrt{D} > 0$, and that the fundamental unit $\epsilon$ is normalized to satisfy $\epsilon > 1$.} Then
\[ \text{$\Oo$ is half-factorial} \Longleftrightarrow \epsilon~\text{generates}~(\Oo_K/f\Oo_K)^{\times}/\langle \text{images of integers prime to $f$}\rangle. \] Finally, if $f=2p$, where $2$ and $p$ are distinct primes inert in $K$, then 
\[ \text{$\Oo$ is half-factorial} \Longleftrightarrow 3\nmid p+1, \text{ and both $\Oo_2$, $\Oo_p$ are half-factorial}.  \]
\end{thmB} 

From the algebraic standpoint, there is no room for improvement in Theorem B; the stated conditions are both necessary and sufficient. But the statistical question of how often half-factorial orders appear ``in the wild'' is far from settled. For this, one must study how often the conditions of Theorem B are satisfied; this is an analytic problem rather than a purely algebraic one. In \cite{coykendall01}, Coykendall posed two conjectures in this direction:
\begin{enumerate}
\item[(A)] Varying both the quadratic field $K$ and the natural number $f$, one encounters infinitely many half-factorial quadratic orders $\Oo_f$. 
\item[(B)] Fix the field $K = \Q(\sqrt{2})$. Then $\Oo_f$ is a half-factorial domain for infinitely many $f \in \N$.
\end{enumerate}

Coykendall's conjectures were recently studied by the second author (P.P.) in \cite{pollack24}. The weaker Conjecture A is proved in full, while the stronger Conjecture B is demonstrated under the assumption of the \textsf{Generalized Riemann Hypothesis} (GRH).\footnote{Throughout, by GRH we mean the assertion that all nontrivial zeros of all number field zeta functions lie on the line $\Re(s)=\frac{1}{2}$.} The methods of \cite{pollack24} were extended in \cite{pollack25} to obtain various related results. Here is a sample.
\begin{thmC} Let $\mathcal{E} = \{1, \frac{3}{2}, 2, \frac{5}{2}, 3, \frac{7}{2}, \dots\} \cup \{\infty\}$.
\begin{enumerate}[label=\normalfont(\roman*)]
\item If $K$ is any quadratic order, then $\rho(\Oo)\in \mathcal{E}$.
\item Under GRH, each element of $\mathcal{E}$ is the elasticity of infinitely many many distinct orders in $\Q(\sqrt{2})$.
\end{enumerate}\end{thmC}

In both \cite{pollack24, pollack25}, the goal is to realize a prescribed elasticity, and the conductors $f$ are constructed accordingly. It is natural to ask what elasticities one sees if instead of constructing $f$ towards a predetermined end, one samples $f$ ``at random'' and records the results. This question was recently considered by the present authors in \cite{FP}.

It requires some care to settle on the right notion of ``at random.'' Let $K$ be a quadratic field, and let $\Oo$ be the order of conductor $f$ in $K$. It follows from a general theorem of Halter-Koch (see \cite[Corollary 4]{HK95}) that $\rho(\Oo) = \infty$ if and only if $f$ is divisible by a prime $p$ that splits (completely) in $K$. For a given $K$, the Chebotarev density theorem guarantees that asymptotically half of all primes $p$ split in $K$. It follows that asymptotically 100\% of natural numbers $f$ have at least one split prime factor, and so $\rho(\Oo_f) = \infty$ asymptotically 100\% of the time. Thus, if we fix a quadratic field $K$, the natural problem is to study the distribution of $\rho(\Oo_f)$ with $f$ sampled uniformly from (only) \textsf{split-free} integers, meaning integers possessing no split prime factors. This is precisely what we do in \cite{FP}. 

\begin{thmD} Let $K$ be a fixed imaginary quadratic field. Let $\Delta_K$ denote the discriminant of $K$. For almost all split-free numbers $f$,
\begin{equation*}
\rho(\Oo_f) = f/(\log{f})^{\frac{1}{2}\log_3{f} + \frac{1}{2}C_K+O((\log_4f)^3/\log_3f)},
\end{equation*}
where $\log_k$ denotes the $k$th iterate of the natural logarithm, and
\[C_K\colonequals\sum_{p>2}\frac{\log p}{(p-1)^2}-1-\frac{\sgn(\Delta_K)1_{\Dd}(\Delta_K)|\Delta_K|\log\Rad(|\Delta_K|)}{\varphi(|\Delta_K|)^2}.\]
Here $\Rad(|\Delta_K|)$ is the product of the distinct primes dividing $\Delta_K$, and $1_{\Dd}$ is the characteristic function of the set of \textsf{prime discriminants}\[
\Dd\colonequals\{-4,\pm8\}\cup\left\{(-1)^{\frac{p-1}{2}}p\colon p>2\text{~is prime}\right\}.   
\]
\end{thmD}

\begin{thmE}[conditional on GRH] Let $K$ be a fixed real quadratic field. For almost all split-free numbers $f$,
\[ \rho(\Oo_f) = (\log{f})^{\frac{1}{2}+O(1/\log_4{f})}. \]
\end{thmE}

In these theorems, \textsf{almost all} means that the proportion of split-free numbers up to $x$ for which the estimate fails tends to $0$, as $x\to\infty$. We remind the reader that $A=O(B)$ means $|A| \le C|B|$ for some \textsf{implied constant} $C$. In these results, as well as all the theorems appearing below, implied constants may depend on the field $K$.

In this paper, we continue our investigations into the distribution of elasticities of orders belonging to a fixed quadratic field. However, instead of looking for the typical size of $\rho(\Oo_f)$, we inquire into the extremal behavior. How large and how small can $\rho(\Oo_f)$ get, in terms of $f$? 

The following are our principal results. Below, $A \ll B$ is synonymous with $A=O(B)$. That is, $|A|\le C|B|$, for a constant $C$ that may depend on $K$. The notation ``$A\gg B$'' indicates that  $B\ll A$.

\begin{thm}[Maximal order, imaginary case]\label{thm:maxorderimaginary} Let $K$ be a fixed imaginary quadratic field. Then $\rho(\Oo_f) \ll f$ for all split-free numbers $f$. Conversely, if $p$ is a prime inert or ramified in $K$, then $\rho(\Oo_p) \gg p$.
\end{thm}

\begin{thm}[Minimal order, imaginary case]\label{thm:minorderimaginary} There are universal positive constants $c_1$, $c_2$, and $f_0$ for which the following holds: Let $K$ be a fixed imaginary quadratic field. For all split-free numbers $f > f_0$, we have 
\begin{equation}\label{eq:lowerboundMOI} \rho(\Oo_f) > (\log{f})^{c_1 \log\log\log{f}}. \end{equation}
On the other hand, there is a sequence of split-free numbers $f$ tending to infinity along which
\begin{equation}\label{eq:upperboundMOI} \rho(\Oo_f) < (\log{f})^{c_2 \log\log\log{f}}.  \end{equation}
\end{thm}

\begin{thm}[Maximal order, real case]\label{thm:maxorderreal} Let $K$ be a real quadratic field. Then \[ \rho(\Oo_f) \ll f/\log{f} \] for all split-free $f>1$. In the opposite direction, GRH implies that for every $\epsilon > 0$, there are infinitely many primes $p$, inert in $K$, with \begin{equation}\label{eq:rholower} \rho(\Oo_p) > p^{\frac{1}{4}-\epsilon}. \end{equation}
\end{thm}

\begin{thm}[Minimal order, real case]\label{thm:minorderreal} Assume GRH. For every real quadratic field $K$, there is a constant $C_K$ with the property that $\rho(\Oo_f) < C_K$ for infinitely many split-free numbers $f$.
\end{thm}

The reader will have noticed the large gap between the (rigorous) upper bound in Theorem \ref{thm:maxorderreal} and the (GRH-conditional) lower bound. As we indicate in a remark following the proof, we believe the upper bound is sharp; in fact, we conjecture that there are infinitely many primes $p$, inert in $K$, with $\rho(\Oo_p) \gg p/\log{p}$. 

\section{The class group, its principal part, and the pre-class group}\label{sec:prelims}
In this section we set up machinery from \cite{FP} that will play an important role in our proofs.

Let $K$ be a quadratic field. For each natural number $f$, we let $I_K(f)$ denote the group of fractional ideals of $K$ generated by the nonzero ideals of $\Oo_K$ comaximal with $f\Oo_K$. We write $P_{K,\Z}(f)$ for the subgroup of $I_K(f)$ generated by principal ideals $\alpha \Z_K$, where $\alpha \in \Oo_K$ satisfies $\alpha \equiv a\pmod{f\Oo_K}$ for an integer $a$ coprime to $f$. The class group $\Cl(\Oo_f)$ is defined as the quotient $I_K(f)/P_{K,\Z}(f)$. Note that when $f=1$, our definition recovers the usual definition of the class group of $\Oo_K=\Oo_1$. 

The next result, which appears as \cite[Lemma 2.1]{FP}, is a variant of Theorem A for quadratic orders. 

\begin{prop}\label{prop:elasticityestimate} Let $K$ be a quadratic field. For each split-free $f\in \N$, 
\[ \frac{1}{2} \Dav\Cl(\Oo_f) \le \rho(\Oo_f)\le \max\left\{1,\frac{1}{2}\Dav \Cl(\Oo_f) + \frac{3}{2}\Omega(f)\right\}.  \]
\end{prop}
\noindent (Here $\Omega(\cdot)$ counts the total number of prime factors, with multiplicity; e.g., $\Omega(-12) = \Omega(30) = 3$.) Proposition \ref{prop:elasticityestimate} is not as precise as Theorem A; except in the case $f=1$, it does not determine the exact value of $\rho(\Oo_f)$. Nevertheless, it will suffice to obtain our statistical results.

In order to access $\Dav \Cl(\Oo_f)$, we find it helpful to ``pull apart'' the class group. We define the \textsf{principal part of the class group}, denoted $\PrinCl(\Oo_f)$, by
\[ \PrinCl(\Oo_f):= (\Oo_K/f\Oo_K)^{\times}/\langle \text{images of integers prime to $f$, units of $\Oo_K$}\rangle. \] Write $P_K$ for $P_{K,\Z}(1)$ (the group of principal fractional ideals). The exact sequence
\[ (\Z/f\Z)^{\times} \times \Oo_K^{\times} \stackrel{\mu}{\longrightarrow} (\Oo_K/f\Oo_K)^{\times} \stackrel{\iota}{\longrightarrow} (I_K(f) \cap P_K)/P_{K,\Z}(f) \longrightarrow 1 \]
allows us to identify $\PrinCl(\Oo_f)$ with the subgroup  $(I_K(f) \cap P_K)/P_{K,\Z}(f)$ of $\Cl(\Oo_f)$; here $\mu$ and $\iota$ are the maps defined by 
\[ \mu((a\bmod{f}, \eta)):= a\eta\bmod{f\Oo_K}\quad\text{and}\quad \iota(\alpha\bmod{f\Oo_K}) = [\alpha \Oo_K]. \]
This identification explains the term ``principal part of the class group.'' 

With the obvious maps, there is a short exact sequence
\begin{equation*} 1 \longrightarrow \conggto{(I_K(f)\cap P_K)/P_{K,\Z}(f)}{\PrinCl(\Oo_f)} \longrightarrow \equalto{I_K(f)/P_{K,\Z}(f)}{\Cl(\Oo_f)} \longrightarrow \equalto{I_K/P_K}{\Cl(\Oo_K)} \longrightarrow 1. \end{equation*}
(Exactness at the last position is not obvious; for this one uses that each ideal class in $\Oo_K$ has a representative comaximal with $f\Oo_K$.) Thus, we may view $\PrinCl(\Oo_f)$ as a subgroup of $\Cl(\Oo_f)$ for which
\begin{equation}\label{eq:indexprin} [\Cl(\Oo_f): \PrinCl(\Oo_f)] = \#\Cl(\Oo_K) = h_K. \end{equation}

To apply Proposition \ref{prop:elasticityestimate} requires knowledge of $\Dav \Cl(\Oo_f)$. Equation \eqref{eq:indexprin} is helpful here; it implies that
\begin{equation}\label{eq:ratio0} \Dav \PrinCl(\Oo_f) \le \Dav \Cl(\Oo_f) \le h_K \Dav \PrinCl(\Oo_f). \end{equation}
(See \cite[Lemma 2.4]{FP}.) In our results, an ambiguity up to a factor depending on $K$ is acceptable, and \eqref{eq:ratio0} allows us to always work with $\PrinCl(\Oo_f)$ rather than $\Cl(\Oo_f)$.  

We will understand $\PrinCl(\Oo_f)$ by viewing it as arising from a two-stage construction. First, we quotient $(\Oo_K/f\Oo_K)^{\times}$ (only) by the image of the integers prime to $f$; we call this is the \textsf{pre-class group}. That is,
\[ \PreCl(\Oo_f) := (\Oo_K/f\Oo_K)^{\times}/\langle \text{images of integers prime to $f$}\rangle. \] (The observant reader will have noticed that this group appeared already in Theorem B.)  Second, we take the quotient of $\PreCl(\Oo_f)$ by the image $\Uu_f$ of $\Oo_K^{\times}$ in $\PreCl(\Oo_f)$; then $\PreCl(\Oo_f)/\Uu_f \cong \PrinCl(\Oo_f)$. 

It is productive to view these pre-class groups $\PreCl(\Oo_f)$ as close cousins of the unit groups $(\Z/m\Z)^{\times}$; this principle is exploited heavily in \cite{FP} (and implicitly in \cite{pollack24, pollack25}). Let $\psi(f):=\#\PreCl(\Oo_f)$. By the Chinese remainder theorem,
\begin{equation}\label{eq:preCRT} \PreCl(\Oo_f) \cong \prod_{p^k\parallel f} \PreCl(\Oo_{p^k}). \end{equation}
Thus, $\psi$ is a multiplicative function of $f$.  Furthermore, since the $\varphi(p^k)$ integers in $[1,p^k]$ that are prime to $p^k$ have distinct images in $\Oo_K/p^k\Oo_K$,
\[ \psi(p^k) = \frac{1}{\varphi(p^k)} (\#\Oo_K/p^k\Oo_K)^{\times} = \frac{1}{p^{k-1}(p-1)} (\#\Oo_K/p^k\Oo_K)^{\times}. \]
Now $\#(\Oo_K/p^k\Oo_K)^{\times} = N(p^k\Oo_K) \prod_{P \mid p}\left(1-\frac{1}{N(P)}\right) = p^{2k} \prod_{P \mid p}\left(1-\frac{1}{N(P)}\right)$. Plugging the result of this formula into the last display, considering separately the cases when $p$ is ramified, split, or inert, we find after a short calculation that
\[ \psi(p^k) = p^k - \leg{\Delta_K}{p} p^{k-1}= p^k \left(1-\leg{\Delta_K}{p}\frac{1}{p}\right), \]
where $\leg{\cdot}{\cdot}$ is the Kronecker symbol.
Hence,
\[ \psi(f) = f \prod_{p\mid f} \left(1-\leg{\Delta_K}{p}\frac{1}{p}\right), \]
in close analogy with the familiar formula $\varphi(f) = f\prod_{p\mid f}(1-\frac1p)$ for Euler's $\varphi$-function.

The groups $\PreCl(\Oo_{f})$ mirror $(\Z/m\Z)^{\times}$ not only in their size but also in their structure. For instance, we have the following result.

\begin{lem}\label{lem:usuallycyclic} Let $K$ be a quadratic field, and let $p$ be a rational prime with $p>3$. Then $\Cl(\Oo_{p^k})$ is cyclic for each natural number $k$. 
\end{lem}

Lemma \ref{lem:usuallycyclic} is an analogue of Gauss's classical theorem that $(\Z/p^k\Z)^{\times}$ is cyclic for each prime power $p^k$ with $p>2$.

\begin{proof}[Proof of Lemma \ref{lem:usuallycyclic}] In the case where the prime $p>3$ is ramified or inert in $K$, this result appears as Lemma 2.5 of \cite{FP}, where it is deduced from related structure theorems of Halter-Koch \cite{HK72}. Suppose now that $p$ is split, say $p \Oo_K = P_1 P_2$ with the $P_i$ distinct maximal ideals of $\Oo_K$. Then $\Oo_K/p^k\Oo_K \cong \Oo_K/P_1^k \times \Oo_K/P_2^k$, and each $\Oo_K/P_i^k \cong \Z/p^k\Z$. Furthermore, 
\begin{align*} \PreCl(\Oo_{p^k}) &= (\Oo_K/p^k\Oo_K)^{\times}/\langle \text{images of integers prime to $p^k$}\rangle\\
&\cong \frac{(\Z/p^k\Z)^{\times} \times (\Z/p^k\Z)^{\times}}{(\Z/p^k\Z)^{\times}} \\
&\cong (\Z/p^k\Z)^{\times};
\end{align*}
here the $(\Z/p^k\Z)^{\times}$ in the ``denominator'' of the second line is viewed as a  subgroup of the ``numerator'' via the diagonal embedding. By the theorem of Gauss quoted above, $(\Z/p^k\Z)^{\times}$ is cyclic, and hence $\PreCl(\Oo_{p^k})$ is as well.
\end{proof}

It is clear from the isomorphism $(\Z/m\Z)^{\times} \cong \prod_{p^k\parallel m} (\Z/p^k\Z)^{\times}$ that the exponent $\Exp{(\Z/m\Z)^{\times}}$ is a divisor of $\lcm\{\varphi(p^k): p^k \parallel m\}$, and Gauss showed in his \emph{Disquisitiones} that the corresponding quotient is always $1$ or $2$. Our proof of Theorem \ref{thm:minorderimaginary} requires the following variant for pre-class groups. Put
\begin{equation}\label{eq:Lfdef} L(f):= \lcm\{\psi(p^k): p^k\parallel f\}. \end{equation}

\begin{prop}\label{prop:Lfdivisibility} For all natural numbers $f$,
\[ \Exp \PreCl(\Oo_f) \mid L(f) \mid 12  \Exp \PreCl(\Oo_f). \]
\end{prop}

It follows from Proposition \ref{prop:Lfdivisibility} that $L(f)/\Exp \PreCl(\Oo_f) \in \{1,2,3,4,6,12\}$. 

\begin{proof} The first divisibility is clear from the isomorphism \eqref{eq:preCRT}, keeping in mind that $\#\PreCl(\Oo_{p^k}) = \psi(p^k)$. To prove the second, it suffices to show that 
\begin{equation*} \psi(p^k) \mid 12 \Exp \PreCl(\Oo_{p^k}) \quad\text{for each prime power $p^k$}, \end{equation*} for then 
\begin{align*} L(f) \mid \lcm \{12 \Exp \PreCl(\Oo_{p^k}): p^k \parallel f\} &=  12 \lcm \{\Exp \PreCl(\Oo_{p^k}): p^k \parallel f\} \\ &= 12 \Exp \PreCl(\Oo_f),\end{align*}
as desired.

When $p>3$, we have $\psi(p^k) = \Exp \PreCl(\Oo_{p^k})$ (Lemma \ref{lem:usuallycyclic}). So we may assume $p=2$ or $p=3$. 

Suppose first that $p=2$. Since $\psi(2^k) = (2-\leg{\Delta_K}{2}) 2^{k-1} \mid 12$ when $k \in \{1,2\}$, we may assume that $k\ge 3$. Let $\alpha = 1+2\sqrt{D}$, where $D$ is the squarefree integer with $K=\Q(\sqrt{D})$. A straightforward induction shows that, for each integer $j\ge 0$,
\[ \alpha^{2^j} = u_j + v_j\sqrt{D} \] 
for integers $u_j, v_j$ with $u_j$ odd and $2^{j+1}\parallel v_j$. In particular, $\alpha^{2^{k-1}} \in \Z[2^{k}\sqrt{D}] \subset \Oo_{2^k}$. Since $\alpha\Oo_K$ and $2\Oo_K$ are comaximal, we conclude that $\alpha$ represents an element of $\PreCl(\Oo_{2^k})$ of order $e$ (say) with \[ e\mid 2^{k-1}.\] On the other hand,
\[ 2 \alpha^{2^{k-3}} = 2u_{k-3} + 2 v_{k-3}\sqrt{D} \notin \Z[2^k\sqrt{D}]. \]
As $2\Oo_{2^k} \subset \Z[2^k\sqrt{D}]$, it must be that $\alpha^{2^{k-3}} \notin \Oo_{2^k}$. Thus, \[ e\nmid 2^{k-3}. \]
Combining the last two displays, we find that $e=2^{k-2}$ or $e=2^{k-1}$. Hence,
\[ \psi(2^k) = (2-\tlegendre{\Delta_K}{2}) 2^{k-1} \mid 12\cdot 2^{k-2} \mid 12e \mid 12 \Exp \PreCl(\Oo_{2^k}). \]

The proof when $p=3$ is similar, and we condense the details. Since $\psi(3^k) = (3-\leg{\Delta_K}{3}) 3^{k-1}\mid 12$ when $k = 1$, we can assume that $k \ge 2$. In this case, one takes $\alpha' = 1+3\sqrt{D}$ and shows by induction that for each integer $j\ge 0$,
\[ \alpha'^{3^j} = {u'_j} + {v'_j}\sqrt{D} \] with integers $u'_j, v'_j$ for which $3\nmid u'_j$, $3^{j+1}\parallel v'_j$. Then $\alpha'^{3^{k-1}} \in \Oo_{3^k}$ while $\alpha'^{3^{k-2}} \notin \Oo_{3^k}$. It follows that $\alpha'$ represents an element in $\PreCl(\Oo_{3^k})$ of exact order $3^{k-1}$. Hence,
\[ \psi(3^k) = (3-\tlegendre{\Delta_K}{3}) 3^{k-1} \mid 12 \cdot 3^{k-1} \mid 12 \Exp \PreCl(\Oo_{3^k}), \]
as desired.\end{proof}

We saw above that $\PrinCl(\Oo_f) \cong \PreCl(\Oo_f)/\Uu_f$, where $\Uu_f$ is the image of $\Oo_K^{\times}$ inside $\PreCl(\Oo_f)$. Hence, setting
\[ \ell(f):= \#\Uu_f, \]
we have that
\[ \#\PrinCl(\Oo_f) = \frac{\psi(f)}{\ell(f)}. \]
We finish this section by recording two observations about $\Uu_f$ that will be relevant in the sequel:
\begin{enumerate}
\item[(i)] When $K$ is imaginary quadratic, $\ell(f) \le \#\Oo_K^{\times} \le 6$. Moreover, $\Uu_f$ is cyclic, generated by (the image of) a primitive $\#\Oo_K^{\times}$-th root of unity.
\item[(ii)] When $K$ is real, $\Oo_K^{\times}$ is generated by $\pm 1$ and the fundamental $\epsilon$. Since $-1$ projects to the identity in $\PreCl(\Oo_f)$, it follows that $\Uu_f$ is cyclic in this case as well, generated by the image of $\epsilon$ in $\PreCl(\Oo_f)$. Furthermore,
\begin{equation*} \ell(f) \text{ is the least integer $\ell$ with $\epsilon^{\ell} \in \Oo_f$}. \end{equation*}
\end{enumerate}

\section{Maximal order in the imaginary case: Proof of Theorem \ref{thm:maxorderimaginary}}\label{sec:maxorderimaginary}
\subsection{Upper bound} We start by showing that $\rho(\Oo_f) \ll f$ for each split-free $f$. We may (and will) assume that $f > 1$. By Proposition \ref{prop:elasticityestimate}, for each split-free number $f>1$,
\begin{equation*} \rho(\Oo_f) \le \frac{1}{2}\Dav \Cl(\Oo_f) + \frac{3}{2}\Omega(f). \end{equation*}
Clearly, 
\begin{equation*} 2^{\Omega(f)} = \prod_{p^k\parallel f} 2^k \le \prod_{p^k\parallel f}p^k  = f, \quad\text{so that}\quad \Omega(f) \ll \log{f}, \end{equation*} and $\log{f} \ll f$. Furthermore, we have from \eqref{eq:ratio0} that
\begin{equation*}\Dav \Cl(\Oo_f) \le h_K \Dav \PrinCl(\Oo_f). \end{equation*}
It therefore suffices to show that $\Dav \PrinCl(\Oo_f) \ll f$. This is an immediate consequence of the following result, which will also be needed later for real $K$.

\begin{prop}\label{prop:davprinupper} Let $K$ be any quadratic field. For every natural number $f$, 
\begin{equation}\label{eq:davprinupper} \Dav \PrinCl(\Oo_f) \ll \frac{f}{\ell(f)}.\end{equation}
\end{prop}

The main ingredient in the proof of Proposition \ref{prop:davprinupper} is an elegant  result of Van Emde Boas and Kruyswijk \cite{EBK67} which bounds from above the Davenport constant of a finite abelian group $G$ in terms of $\Exp G$.
\begin{prop}\label{prop:davexp} Let $G$ be a finite abelian group, and put 
\[ r(G) = \frac{\#G}{\Exp G}.\]
Then 
\[ \Dav G \le \#G \left(\frac{1+\log r(G)}{r(G)}\right).\]
\end{prop}

The function $t\mapsto \frac{1+\log{t}}{t}$ is decreasing on the domain $t\ge 1$, with maximum value $1$. Thus, $D(G) \le \#G$ always, with strict inequality whenever $r(G) > 1$ (that is, whenever $G$ is not cyclic).

In order for Proposition \ref{prop:davexp} to be of use in the proof of \eqref{eq:davprinupper}, we need a handle on $r(\PrinCl(\Oo_f))$. This will be given to us by the next lemma. Below, we write $\rk_2 G$ for the $2$-rank of $G$.

\begin{lem}\label{lem:rexp} For every finite abelian group $G$,
\[\exp G \mid 2^{1-\rk_2 G} \#G. \]
Hence, 
\[ r(G) \ge 2^{\rk_2{G} - 1}. \]
\end{lem}
\begin{proof} The second claim is immediate from the first, so we focus on that.

If $\rk_2 G \le 1$, the asserted divisibility is trivial, so we assume $\rk_2 G > 1$. We may also assume $G = A \oplus \Z/2^{r_1}\Z \oplus \Z/2^{r_2}\Z \oplus \dots \oplus \Z/2^{r_l}\Z$, where $l = \rk_2 G$, the $r_i$ satisfy $1 \le r_1 \le r_2 \le \dots \le r_l$, and $A$ has odd order. Then
\[ \Exp G \mid 2^{r_{l}} \#A = 2^{1-l} \cdot 2^{(l-1) + r_{l}} \#A \mid 2^{1-l} \cdot 2^{r_1 + r_2 + \dots + r_{l-1} +r_{\ell}}\#A = 2^{1-l} \#G. \qedhere\]
\end{proof}

\begin{proof}[Proof of Proposition \ref{prop:davprinupper}] Let $w$ denote the number of (distinct) odd primes $p$ dividing $f$ which are inert in $K$. We will bound $\rk_2 \PrinCl(\Oo_f)$ below in terms of $w$ and then apply Lemma \ref{lem:rexp}.

If $p^k\parallel f$ with $p$ odd and inert in $K$, then $\#\PreCl(\Oo_{p^k}) =\psi(p^k) = p^{k} + p^{k-1}$ is even. Hence, $\rk_2 \PreCl(\Oo_{p^k}) \ge 1$, and 
\begin{equation*} \rk_2 \PreCl(\Oo_f) = \sum_{p^k\parallel f} \rk_2 \PreCl(\Oo_{p^k}) \ge w. \end{equation*}
As $\Uu_f$ is cyclic, $\rk_2{\Uu_f} \le 1$, and 
\[ \rk_2 \PrinCl(\Oo_f) = \rk_2 \frac{\PreCl(\Oo_f)}{\Uu_f} \ge \rk_2 \PreCl(\Oo_f) - \rk_2 \Uu_f \ge w-1.\]
Thus, by Lemma \ref{lem:rexp},
\begin{equation*} r:= r(\PrinCl(\Oo_f)) \ge 2^{w-2}. \end{equation*}

Suppose that $w\ge 2$. Then $r \ge 2^{w-2}\ge 1$, and Proposition \ref{prop:davexp} yields
\[ \Dav\PrinCl(\Oo_f) \le \frac{\psi(f)}{\ell(f)} \frac{1+\log{r}}{r} \le \frac{\psi(f)}{\ell(f)} \frac{1+\log{2^{w-2}}}{2^{w-2}} < \frac{4w}{2^w} \cdot \frac{\psi(f)}{\ell(f)}. \]
As $f$ is split-free, we have $\omega(f) \le w+1+\omega(\Delta_K)$. (Per the usual convention in analytic number theory, we write $\omega(\cdot)$ for the count of distinct prime factors.) Therefore,
\[ \psi(f) \le f\prod_{p\mid f} \left(1+\frac1p\right) \le (3/2)^{\omega(f)} f \le (3/2)^{w+1+\omega(\Delta_K)} f. \] and
\[ \Dav\PrinCl(\Oo_f) < \frac{4w}{2^w} \frac{\psi(f)}{\ell(f)} \le (3/2)^{1+\omega(\Delta_K)} (4w (3/4)^w)\frac{f}{\ell(f)}. \]
In this last expression, the coefficient of $f/\ell(f)$ is bounded by a constant depending only $K$, establishing \eqref{eq:davprinupper} when  $w\ge 2$.

If $w\le 1$, the argument is easier. In this case, $\omega(f) \le 2+ \omega(\Delta_K)$, and
\begin{align*} \Dav \PrinCl(\Oo_f) &\le \#\PrinCl(\Oo_f) = \frac{\psi(f)}{\ell(f)}\\
&=\frac{f}{\ell(f)} \prod_{p \mid f} \left(1-\frac{1}{p}\leg{\Delta_K}{p}\right)\le (3/2)^{\omega(f)} \frac{f}{\ell(f)} \ll \frac{f}{\ell(f)},
\end{align*}
as desired. 
\end{proof}

\subsection{Lower bound} We turn now to the (much simpler) proof that $\rho(\Oo_p) \gg p$ for all $p$ that do not split completely in $K$. If $p =2$ or $p=3$, then $\rho(\Oo_p) \ge 1 \ge \frac{1}{3}p$, so we suppose that $p>3$. By Proposition \ref{prop:elasticityestimate} and \eqref{eq:ratio0}, 
\[ \rho(\Oo_p) \ge \frac{1}{2} \Dav \PrinCl(\Oo_p). \]
The group $\PrinCl(\Oo_p)$ is cyclic, as a quotient of the cyclic group $\PreCl(\Oo_p)$ (see Lemma \ref{lem:usuallycyclic}). Whenever $G$ is cyclic, $\Dav G = \#G$ (see, e.g., Lemma 1.4.9 on p.\ 27 of \cite{GHK06}.) Thus, $\Dav \PrinCl(\Oo_p) = \#\PrinCl(\Oo_p) = \psi(p)/\ell(p)$, 
and we conclude that
\begin{equation}\label{eq:rhoplower} \rho(\Oo_p) \ge \frac{\psi(p)}{2\ell(p)}. \end{equation}

Up to this point all of our reasoning has been valid whether $K$ is real or imaginary. But if we assume $K$ is imaginary, then $\ell(p)  \le \#\Oo_K^{\times} \le 6$, leading to the conclusion that
\[ \rho(\Oo_p) \ge \frac{\psi(p)}{12} = \frac{p-\leg{\Delta_K}{p}}{12} \ge \frac{p}{12} \]
for all non-split $p>3$. This establishes the lower bound claimed in Theorem \ref{thm:maxorderimaginary}. In fact, we have shown that the implied constant there can be taken as $\frac{1}{12}$, independently of $K$.

\section{Minimal order in the imaginary case: Proof of Theorem \ref{thm:minorderimaginary}}

Our proof of Theorem \ref{thm:minorderimaginary} is based on arguments of Erd\H{o}s, Pomerance, and Schmutz \cite{EPS91} used to estimate the minimal order of \textsf{Carmichael's lambda function} $\lambda(m)\colonequals\#\Exp (\Z/m\Z)^{\times}$. According to the first half of Theorem 1 of \cite{EPS91}, we have \[ \lambda(m) > (\log{m})^{(\frac{1}{\log{2}} + o(1))\log\log\log{m}} \]
whenever $m\to\infty$. (As usual in analytic number theory, $o(1)$ denotes a quantity tending to $0$.) The second half of that theorem asserts the existence of a strictly increasing sequence of natural numbers $\{m_i\}_{i\ge 1}$ such that
\[ \lambda(m_i) < (\log{m_i})^{c_0 \log\log\log{m_i}} \]
for all $i\ge 1$.

We need an $L(f)$-analogue of the quoted results, where $L(f)$ is the function defined in \eqref{eq:Lfdef}.

\begin{prop}\label{prop:minorderL(f)}
Let $K$ be a quadratic field. For any sequence of natural numbers $f\to\infty$,
\[L(f) > (\log f)^{(\frac{1}{\log{2}}+o(1))\log\log\log f}, \]
uniformly in $K$. On the other hand, for some universal constant $c_0>0$, and each fixed choice of $K$, there is a strictly increasing sequence $\{f_i\}_{i\ge1}$ of positive integers, each of which is a product of distinct primes inert in $K$, such that
\[L(f_i)<(\log f_i)^{c_0\log\log\log f_i}\]
for all $i\ge1$. 
\end{prop}

Taking Proposition \ref{prop:minorderL(f)} as shown, we can quickly conclude the proof of Theorem \ref{thm:minorderimaginary}. We start with the lower bound \eqref{eq:lowerboundMOI}. From Proposition \ref{prop:elasticityestimate} and \eqref{eq:ratio0}, $\rho(\Oo_f) \ge \frac{1}{2}\Dav \PrinCl(\Oo_f)$. As $\PrinCl(\Oo_f)$ is the quotient of $\PreCl(\Oo_f)$ by the subgroup $\Uu_f$ of order $\#\Uu_f \le \#\Oo_K^{\times} \le 6$, we have $\Dav \PrinCl(\Oo_f) \ge \frac{1}{6} \Dav \PreCl(\Oo_f)$ (see \cite[Lemma 2.4]{FP}). Therefore, 
\[ \rho(\Oo_f) \ge \frac{1}{12} \Dav \PreCl(\Oo_f) \ge \frac{1}{12} \Exp \PreCl(\Oo_f) \ge \frac{1}{12^2} L(f), \]
invoking Proposition \ref{prop:Lfdivisibility} at the last step. The lower bound of Proposition \ref{prop:minorderL(f)} thus implies that \eqref{eq:lowerboundMOI} holds for any constant $c_1 < 1/\log{2}$. 
 
The upper bound \eqref{eq:upperboundMOI} is slightly more intricate. From Proposition \ref{prop:elasticityestimate}, \eqref{eq:ratio0}, and the bound $\Omega(f)\ll \log{f}$, we have \begin{equation}\label{eq:upperMOI1}\rho(\Oo_f) \ll \Dav \PrinCl(\Oo_f) + \log{f}\end{equation} for all split-free $f$. Applying the Van Emde Boas--Kruyswijk exponent bound (Proposition \ref{prop:davexp}), 
\begin{align*} \Dav \PrinCl(\Oo_f) &\le \Exp \PrinCl(\Oo_f) \left(1 + \log \frac{\#\PrinCl(\Oo_f)}{\Exp \PrinCl(\Oo_f)}\right) \\
&\le \Exp \PrinCl(\Oo_f) (1 + \log{\psi(f)}). \end{align*}
Now $\psi(f) \le (3/2)^{\omega(f)} f$, so that $$ \log \psi(f) \le \omega(f) \log\frac{3}{2} + \log{f} \le \Omega(f) \log\frac{3}{2} + \log{f} \ll \log{f}.$$ We conclude that for all split-free $f$,
\begin{align} \Dav \PrinCl(\Oo_f) &\ll \Exp \PrinCl(\Oo_f) (1+\log{f})\notag \\
&\le L(f) (1+\log{f}). \label{eq:upperMOI2}\end{align}
The upper bound half of Theorem \ref{thm:minorderimaginary} now follows from \eqref{eq:upperMOI1}, \eqref{eq:upperMOI2}, and the upper bound result of Proposition \ref{prop:minorderL(f)}. Indeed, fixing any $c_2 > c_0$, we find that \eqref{eq:upperboundMOI} holds for all large enough $f$ (large enough in terms of $K$) from the sequence $\{f_i\}$.

The remainder of this section is devoted to the proof of Proposition \ref{prop:minorderL(f)}. We make little claim to originality here; our arguments are straightforward adaptations of those in \cite{EPS91}.

\subsection{The lower bound in Proposition \ref{prop:minorderL(f)}}

We start by observing that $L(f)\to\infty$ whenever $f\to\infty$. Indeed, if $L(f)\le B$ for some constant $B>0$, then $p^k\le 2B$ for all $p^k\parallel f$, since $L(p^k)\ge \varphi(p^k)\ge p^k/2$ for all prime powers $p^k$. Consequently, there are only finitely many $f$ for which $L(f)\le B$. 

Next, a case-by-case analysis shows that $L$ is at most 4-to-1 on prime powers. To see this, we show that every $m\in\N$ is the $L$-image of at most $4$ prime powers. Since the restriction of $L$ to powers of ramified primes is the identity function, the $L$-preimage of $m$ contains at most one power of a ramified prime. Let $p$ be the least prime for which there exists some $k\in\N$ with $L(p^k)=m$.

\underline{Case I. $p$ is ramified in $K$.}

Let $L(q^l)=m$ for some $q^l\ne p^k$, where $q>p$ is unramified in $K$. Suppose first that $q$ splits in $K$. For this case, we have $p^k=q^{l-1}(q-1)$, which implies that $l=1$ and thus $q=p^k+1$. Suppose now that $q$ is inert in $K$. Then $p^k=q^{l-1}(q+1)$. So $l=1$ and $q=p^k-1$. We conclude that $m$ is the $L$-image of at most $3$ prime powers. 

\underline{Case II. $p$ splits in $K$.}

Again, let $L(q^l)=m$ for some $q^l\ne p^k$ with $q\ge p$. Suppose first that $q$ is ramified in $K$. Then $q>p$, which leads to $p^{k-1}(p-1)\ne q^l$, contradicting $L(p^k)=L(q^l)$. If $q$ splits in $K$, then $p^{k-1}(p-1)=q^{l-1}(q-1)$. Since $q^l\ne p^k$, we have $q>p$. This implies that $l=1$ and $q=p^{k-1}(p-1)+1$. If $q$ is inert in $K$, then $q>p$ and $p^{k-1}(p-1)=q^{l-1}(q+1)$. Once again, we must have $l=1$, from which we deduce that $q=p^{k-1}(p-1)-1$. So as in Case I, $m$ is the $L$-image of at most $3$ prime powers. 

\underline{Case III. $p$ is inert in $K$.}

Suppose that $L(q^l)=m$ for some $q^l\ne p^k$ with $q\ge p$. If $q$ is ramified in $K$, then $q>p$ and $p^{k-1}(p+1)=q^l$. In particular, we have $q\ge p+1$ and $q^l\mid(p+1)$. So we must have $l=1$ and $q=p+1$. Next, if $q$ splits in $K$, then $q>p$ and $p^{k-1}(p+1)=q^{l-1}(q-1)$. If $l>1$, then the same argument shows that $l=2$ and $q=p+1$. Of course, we get $q=p^{k-1}(p+1)+1$ when $l=1$. Finally, if $q$ is inert in $K$, then $p^{k-1}(p+1)=q^{l-1}(q+1)$. Since $q^l\ne p^k$, we have $q>p$. When $l>1$, we get $l=2$ and $q=p+1$ once again. When $l=1$, it is obvious that $q=p^{k-1}(p+1)-1$. In view of the fact that $q=p+1$ can occur for at most one of the three possible cases for $q$, we conclude that $m$ is the $L$-image of at most $4$ prime powers.

This concludes the proof that $L$ is at most 4-to-1 on prime powers.

We can now complete the proof of the lower bound on $L(f)$. Let $f$ be sufficiently large, and suppose that $L(f)=m$. Since $p^k\le 2\varphi(p^k)\le 2L(p^k)$ for all prime powers $p^k$, we have
\[f\le \prod_{\substack{p^k\mid f\\L(p^k)\mid m}}p^k\le \prod_{d\mid m} \prod_{p^k:\,L(p^k)=d} p^k \le \prod_{d\mid m} (2m)^4=(2m)^{4\tau(m)},\]
where $\tau(\cdot)$ is the divisor counting function. Using the inequality $\tau(m)\le 2^{(1+o(1))\log m/\log\log m}$ (see \cite[Theorem 317, p.\ 345]{HW08}), we deduce that
\[f\le\exp\left((4\log(2m))2^{(1+o(1))\log m/\log\log m}\right),\]
from which it follows that 
\[L(f)=m\ge(\log f)^{(1/\log 2+o(1))\log\log\log f}.\]
As the inequality ``$f \le (2m)^{4\tau(m)}$'' holds independently of $K$, the ``$o(1)$'' term appearing in these last two displays decays to $0$ uniformly in $K$.

\subsection{The upper bound in Proposition \ref{prop:minorderL(f)}}
We need two lemmas. The first of these appears as Proposition 8 in \cite{APR83}. For each real number $x>0$ and each pair of integers $a, k$ with $k \in \N$ and $\gcd(a,k)=1$, put
\[ \theta(x;k,a) = \sum_{\substack{p \le x \\ p\equiv a\psmod{k}}} \log{p}. \]

\begin{lem}\label{lem:PNTAP} Let $\epsilon > 0$. There are computable positive constants $\delta$ and $x_0$ such that
\[ \left|\theta(x;k,a) - \frac{x}{\varphi(k)}\right| < \epsilon \frac{x}{\varphi(k)}\]
for all $x\ge x_0$, all $k\in\N\cap[1,x^{\delta}]$ and all $a\in\Z$ with $\gcd(k,a)=1$, except possibly for those $k$ which are divisible by a certain integer $k_0=k_0(x)>(\log x)^{3/2}$. 
\end{lem}

The next result is a variant of \cite[Proposition 10]{APR83}.

\begin{lem}\label{lem:mxconstruction} Let $K$ be a quadratic field. For all $x$ that are sufficiently large in terms of $K$, there is a natural number $M_x \le x^2$ for which
\begin{equation}\label{eq:omega*inert}
\#\{p~\text{inert in}~K\colon (p+1)\mid M_x\}>\exp\left(C\frac{\log x}{\log\log x}\right). 
\end{equation}
Here $C$ is a positive, universal constant.
\end{lem}

\begin{proof} We let $\delta,x_0>0$ be the parameters of Lemma \ref{lem:PNTAP} corresponding to the choice $\epsilon=\frac12$. 

Without loss of generality, we may assume that $x_0$ is sufficiently large (in terms of $K$). For $x\ge x_0$, let $k_1=k_1(x)$ denote the product of all the unramified primes $p\le \frac{1}{2}\delta\log x$.
Then $k_1<x^{\delta}/|\Delta_K|$. Let $P$ be the largest prime factor of $\Delta_K$, and let $p_0\ge P+1$ be an arbitrary common prime factor of $k_0$ and $k_1$, with the convention that $p_0=1$ when $\gcd(k_0,k_1)$ has no prime factors at least $P+1$. Put $k=k_1/p_0$. If $p_0>1$, then it is obvious that $k_0\nmid k\Delta_K$. On the other hand, if $p_0=1$ and $k_0\mid k\Delta_K$, then all of the prime factors of $k_0$ would be at most $P$. Since $k_0>(\log x)^{3/2}$, there would exist some $p\le P$ such that $p^4\mid k_0$. But $\gcd(k,\Delta_K)=1$, $k$ is squarefree, and $\Delta_K=2^rs$, where $r\in\{0,2,3\}$ and $s\in\Z$ is odd and squarefree, so that $p^4\nmid k\Delta_K$, a contradiction. Therefore, we have $k_0\nmid k\Delta_K$ in both cases.

For each $d\mid k$, we have from Lemma \ref{lem:PNTAP} that
\begin{align*}
\sum_{\substack{p\le x\\p\equiv-1\psmod{d}\\p\text{~inert in~}K}}1&=\sum_{\substack{j\in(\Z/|\Delta_K|\Z)^{\times}\\(\Delta_K/j)=-1}}\sum_{\substack{p\le x\\p\equiv a_j\psmod{|\Delta_K|d}}}1\\
&\ge\frac{1}{\log x}\sum_{\substack{j\in(\Z/|\Delta_K|\Z)^{\times}\\(\Delta_K/j)=-1}}\theta(x;|\Delta_K|d,a_j)\\
&>(1-\epsilon)\frac{x}{\varphi(|\Delta_K|d)\log x}\sum_{\substack{j\in(\Z/|\Delta_K|\Z)^{\times}\\(\Delta_K/j)=-1}}1\\
&=(1-\epsilon)\frac{x}{\varphi(|\Delta_K|)\varphi(d)\log x}\cdot\frac{\varphi(|\Delta_K|)}{2}\\
&=\frac{x}{4\varphi(d)\log x},
\end{align*}
where $a_j\in(\Z/|\Delta_K|d\Z)^{\times}$ satisfies $a_j\equiv -1\psmod{d}$ and $a_j\equiv j\psmod{|\Delta_K|}$.

We now use this estimate to lower bound the cardinality of the set $A$ of pairs $(m,p)\in(\N\cap[1,x])^2$ satisfying the congruence $m(p+1)\equiv0\psmod{k}$, where $p$ is inert in $K$. To this end, we define 
\[A_d\colonequals\left\{(m,p)\in A\colon d\mid(p+1)\text{~and~}\gcd(m,k)=k/d\right\}\]
for each $d\mid k$. Then $A$ is the disjoint union of $A_d$'s over $d\mid(p+1)$. Note that the number of $m\le x$ with $\gcd(m,k)=k/d$ is at least $\varphi(d)\lfloor x/k\rfloor$. Hence, we obtain
\[\#A_d>\frac{x}{4\varphi(d)\log x}\cdot\varphi(d)\lfloor x/k\rfloor>\frac{x^2}{5k\log x}\]
and 
\[\#A=\sum_{d\mid k}\#A_d>\frac{x^2}{5k\log x}\cdot\tau(k)=\frac{x^2}{5k\log x}\cdot2^{\omega(k)}>\frac{x^2}{5k\log x}\cdot2^{\frac{(\delta/4)\log x}{\log\log x}},\]
where we have used the prime number theorem to get
\[\omega(k)\ge\pi\left(\frac{\delta}{2}\log x\right)-\omega(\Delta_K)-1>\frac{(\delta/4)\log x}{\log\log x}\]
for sufficiently large $x$ depending on $\Delta_K$. Since $m(p+1)\in\{n\le x^2\colon k\mid n\}$ for every pair $(m,p)\in A$, some $n\le x^2$ with $k\mid n$ must admit at least 
\[\frac{\#A}{x^2/k}>\frac{1}{5\log x}\cdot2^{\frac{(\delta/4)\log x}{\log\log x}}>\exp\left(\frac{\delta}{6}\frac{\log x}{\log\log x}\right)\]
representations of the form $n=m(p+1)$ with $(m,p)\in A$. (We use here that $\frac{1}{4}\log{2} > \frac16$.) The proof of \eqref{eq:omega*inert} is completed by taking $M_x$ to be this $n$.\end{proof}

We are now in a position to establish the upper bound in Proposition \ref{prop:minorderL(f)}. Let $C$ be as in Lemma \ref{lem:mxconstruction}. Put $x_i=(\log i)^{(2/C)\log\log\log i}$ and 
\[g_i=\prod_{\substack{p~\text{inert in~}K\\(p+1)\mid M_{x_i}}}p.\]
By \eqref{eq:omega*inert}, we have
\[g_i\ge\prod_{\substack{p~\text{inert in~}K\\(p+1)\mid M_{x_i}}}2>\exp\left((\log2)\exp\left(C \frac{\log x_i}{\log\log x_i}\right)\right)>i\]
for sufficiently large $i$. Moreover, we have $L(g_i)\mid M_{x_i}$, which implies that
\[L(g_i)\le M_{x_i}\le x_i^2=(\log i)^{(4/C)\log\log\log i}<(\log g_i)^{c_0\log\log\log g_i}\]
for sufficiently large $i$, where $c_0:=4/C$. The asserted upper bound follows by extracting a strictly monotonic subsequence $\{f_i\}_{i\ge1}$ from $\{g_i\}_{i\ge1}$.

\section{Minimal order in the real case: Proof of Theorem \ref{thm:minorderreal}}

Let $K$ be a real quadratic field. We remind the reader that $\epsilon$ denotes the fundamental unit of $K$, and we set
\[ \delta = \begin{cases}
1 &\text{if $N_{K/\Q}(\epsilon)=1$}, \\
2 &\text{if $N_{K/\Q}(\epsilon)=-1$}.
\end{cases}\]

If $p$ is a prime inert in $K$, then its associated Frobenius element is conjugation on $K$ (the nontrivial element of $\mathrm{Gal}(K/\Q)$). Hence, 
\[ \epsilon^{p+1} \equiv \epsilon^p \epsilon \equiv N_{K/\Q}(\epsilon) \pmod{p\Oo_K}, \]
and
\[ \epsilon^{\delta(p+1)} \equiv N_{K/\Q}(\epsilon)^{\delta} \equiv 1 \pmod{p\Oo_K}. \]
Thus, the order of $\epsilon$ in $(\Oo_K/p\Oo_K)^{\times}$ is a divisor of $\delta(p+1)$. We will base our proof of Theorem \ref{thm:minorderreal} on the following result of Roskam \cite{roskam00}. 
\begin{prop}[conditional on GRH]\label{prop:roskam} There are infinitely many primes $p$, inert in $K$, for which the order of $\epsilon$ in $(\Oo_K/p\Oo_K)^{\times}$ is precisely $\delta(p+1)$.
\end{prop}
\noindent (In fact, Roskam shows that the order is $\delta(p+1)$ not only for infinitely many inert primes, but for a positive proportion of all inert primes. The weaker version here is sufficient for our purposes.) Theorem \ref{thm:minorderreal} is an immediate consequence of Proposition \ref{prop:roskam} in conjunction with the next assertion.

\begin{prop}
If $p$ is a prime inert in $K$ for which $\epsilon$ has order $\delta(p+1)$ in $(\Oo_K/p\Oo_K)^{\times}$, then 
\begin{equation*}\rho(\Oo_p) \le h_K + \frac{3}{2}.\end{equation*}
\end{prop}

\begin{proof} Let $p$ be as in the proposition. Then  $\epsilon^{\ell(p)}\equiv n\pmod{p\Oo_K}$ for some rational integer $n$ prime to $p$. By Fermat's little theorem,
\[ \epsilon^{(p-1)\ell(p)} \equiv 1 \pmod{p\Oo_K}. \]
We are assuming that $\epsilon$ has order $\delta(p+1)$ in $(\Oo_K/p\Oo_K)^{\times}$. Hence, the displayed congruence forces
\[ p+1 \mid \delta(p+1) \mid (p-1)\ell(p). \]
Writing $(p-1)\ell(p) = (p+1)\ell(p) - 2\ell(p)$, we deduce that $p+1 \mid 2\ell(p)$. In particular, \[ \ell(p) \ge \frac{p+1}{2}.\]

By Proposition \ref{prop:elasticityestimate} and \eqref{eq:ratio0},
\[ \rho(\Oo_p) \le \frac{1}{2} \Dav \Cl(\Oo_p) +\frac{3}{2} \le \frac{h_K}{2} \Dav \PrinCl(\Oo_p) + \frac{3}{2} \le \frac{h_K}{2} \#\PrinCl(\Oo_p) + \frac{3}{2}.\]
The proof is completed by observing that $\#\PrinCl(\Oo_p) = \frac{\psi(p)}{\ell(p)} = \frac{p+1}{\ell(p)} \le 2$. 
\end{proof}

\section{Maximal order in the real case: Proof of Theorem \ref{thm:maxorderreal}}\label{sec:maxorderreal}
\subsection{Upper bound} We start by showing that $\rho(\Oo_f) \ll f/\log{f}$ for each $f>1$. By Proposition \ref{prop:elasticityestimate} and \eqref{eq:ratio0}, 
\[ \rho(\Oo_f) \le \frac{h_K}{2} \Dav \PrinCl(\Oo_f) + \frac{3}{2}\Omega(f). \]
As $\Omega(f) \ll \log{f} \ll \frac{f}{\log{f}}$, it suffices to demonstrate that $\Dav \PrinCl(\Oo_f) \ll \frac{f}{\log{f}}$. According to \eqref{eq:davprinupper}, we have $\Dav \PrinCl(\Oo_f) \ll f/\ell(f)$, so we will be done if we prove that \[ \ell(f) \gg \log{f}.\] 

For this we may assume $f$ is sufficiently large (bounded $f$ can be dealt with by adjusting the implied constant). Say $K = \Q(\sqrt{D})$ with $D$ squarefree, and let $\tau_D $ be defined as in \eqref{eq:taudef}, so that $1, f\tau_D$ are a $\Z$-basis for $\Oo_f$. If we express $\epsilon^{\ell(f)} = U + V \tau_D$, then $V$ is a positive integer multiple of $f$. Writing $\epsilon^{\ell(f)} = u + v\sqrt{D}$ (with $u,v\in \Z$), we find that $v \ge \frac{1}{2}V \ge \frac{1}{2}f$. Therefore, using a tilde for conjugation in $K$,
\begin{equation}\label{eq:exponentialgrowth} f \le 2v = \frac{\epsilon^{\ell(f)} - {\tilde{\epsilon}}^{\ell(f)}}{\sqrt{D}} < \frac{\epsilon^{\ell(f)} + 1}{\sqrt{D}}. \end{equation}
Hence (continuing to assume $f$ is large),
\[ \epsilon^{\ell(f)} > f\sqrt{D}-1 > f,\]
and $\ell(f) > \frac{\log{f}}{\log\epsilon} \gg \log{f}$.

\subsection{Lower bound}
Now we turn to the lower bound \eqref{eq:rholower}. Let $p$ be a prime inert in $K$ with $p> 3$. According to \eqref{eq:rhoplower}, 
\[\rho(\Oo_p) \ge \frac{1}{2}\frac{p+1}{\ell(p)}, \]
and so it suffices to produce infinitely many inert 
$p$ with 
\begin{equation}\label{eq:neededratio} \frac{p+1}{\ell(p)} > 2p^{\frac{1}{4}-\epsilon}. \end{equation}
For this we borrow Lemma 6.4 from \cite{FP}.

\begin{lem}[conditional on GRH] Let $K=\Q(\sqrt{D})$ be a real quadratic field, where $D>1$ is squarefree. Let $\eta = \epsilon$ or $\epsilon^2$, according to whether $N_{K/\Q}(\epsilon)=1$ or $-1$, respectively. Let $q$ be an odd prime not dividing $D$, and let $y\ge2$ be a real number. The count of primes $p \le y$ for which 
\begin{equation}\label{eq:pconds} p\text{ is inert in $K$}, \quad p\equiv -1\pmod{q}, \quad\text{ and }\quad \eta^{(p+1)/q}\equiv 1\pmod{p},
\end{equation} is 
\[ \frac{1}{2q(q-1)} \int_{2}^{y}\frac{\mathrm{d}t}{\log{t}} + O(y^{1/2}\log{(qy)}), \]
where the implied constant depends at most on $K$.
\end{lem}

\begin{proof}[Proof that \eqref{eq:neededratio} holds for infinitely many inert $p$] Fix $\epsilon \in (0,1/2)$. For large $y$, we select an arbitrary prime $q$ satisfying
\[ y^{\frac{1}{4}-\frac{1}{2}\epsilon} < q \le y^{\frac{1}{4}-\frac{1}{4}\epsilon}. \]
By the prime number theorem, there will be many choices for $q$ once $y$ is large enough, and none of these will divide $D$. Using $C$ for a constant depending at most on $K$, the count of $p \le y$ satisfying \eqref{eq:pconds} is bounded below by
\begin{align*} \frac{1}{2q(q-1)}\int_{2}^{y} \frac{\mathrm{d}t}{\log{t}} - C y^{1/2}\log{y} &> \frac{1}{2q(q-1) \log{y}}\int_{2}^{y} \mathrm{d}t - C y^{1/2}\log{y} \\
&> \frac{1}{3} \frac{y^{\frac{1}{2}+\frac{1}{2}\epsilon}}{\log{y}} - C y^{1/2}\log{y} 
\\&> y^{\frac12+\frac{1}{4}\epsilon}.
\end{align*}
In particular, for all large $y$ there is at least one such $p$. As $\eta^{(p+1)/q} \equiv 1\pmod{p\Oo_K}$, and $\eta=\epsilon$ or $\epsilon^2$, we must have $\epsilon^{(p+1)/q}\equiv \pm 1\pmod{p\Oo_K}$. Hence, $\epsilon^{(p+1)/q} \in \Oo_p$ and $\ell(p) \mid \frac{p+1}{q}$. Therefore, $q \mid \frac{p+1}{\ell(p)}$, and
\[ \frac{p+1}{\ell(p)} \ge q > y^{\frac{1}{4}-\frac{1}{2}\epsilon} > 2y^{\frac{1}{4}-\epsilon} \ge 2p^{\frac{1}{4}-\epsilon}, \]
verifying \eqref{eq:neededratio}. 

To see that this argument produces infinitely many $p$, note that $p > q > y^{\frac{1}{4}-\frac{1}{2}\epsilon}$, while $y$ can be taken arbitrarily large. \end{proof}

 \begin{rmk} Let $K$ be a real quadratic field. For each $m=0,1,2,\dots$, write $2\epsilon^m = u_m + v_m\sqrt{D}$, where $u_m, v_m \in \Q$. Then  $u_m, v_m \in \N$, and $v_m$ grows exponentially with $m$; in fact, $v_m \sim \epsilon^{m}/\sqrt{D}$, as $m\to\infty$ (compare with \eqref{eq:exponentialgrowth}). 
 
We conjecture the existence of a constant $\delta =\delta_K > 0$ for which the following holds: 
\begin{equation}\label{eq:quadraticmersenne}\text{$v_m$ has an inert prime factor $p > v_m^{\delta}$ for infinitely many  $m\in\N$.}\end{equation}

Suppose $p, m$ are related as in \eqref{eq:quadraticmersenne}, and that $m$ is large. Then $p > 3$, and $p\mid v_m$ implies that $\epsilon^{m} \in \Oo_p$. It follows that \[ \ell(p)\le m \ll \log{v_m} \ll \log{p}, \]  and 
\[ \rho(\Oo_p) \ge \frac{1}{2} \frac{p+1}{\ell(p)} \gg \frac{p}{\log{p}}.\] So if our conjecture \eqref{eq:quadraticmersenne} holds, then the upper bound of Theorem \ref{thm:maxorderreal}(a) is sharp. 

To illustrate, let $K=\Q(\sqrt{2})$. In this case, $\epsilon = 1+\sqrt{2}$, each $v_m$ is even, and $\frac{1}{2}v_m$ is commonly termed the $m$th \textsf{Pell number}. It seems plausible to conjecture that there are infinitely many Pell numbers that are themselves prime; e.g., looking at indices $m<1000$, one finds that $\frac{1}{2}v_m$ is prime for $m=2$, $3$, $5$, $11$, $13$, $29$, $41$, $53$, $59$, $89$, $97$, $101$, $167$, $181$, $191$, $523$, $929$ (this is OEIS sequence \texttt{A096650}). Among these $m$, the prime $\frac{1}{2}v_m$ is inert in $K$ roughly half the time, for $m=3$, $5$, $11$, $13$, $29$, $53$, $59$, $101$, $181$, $523$. We view this as compelling evidence that when $K=\Q(\sqrt{2})$, our hypothesis \eqref{eq:quadraticmersenne} holds for any $\delta \in (0,1)$.
\end{rmk}

\section*{Acknowledgements} P.P. gratefully acknowledges the support of the National Science Foundation under award DMS-2001581.

\bibliographystyle{amsalpha}
\bibliography{elasticity}

\end{document}